\documentclass[twoside,12pt,leqno]{amsart}
\usepackage{pstricks,amssymb,latexsym,verbatim} % verbatim for ``comment''
\usepackage{hyperref}
\hypersetup{citecolor=purple, linkcolor=blue, colorlinks=true}

\theoremstyle{plain}
\newtheorem{theorem}{Theorem}%[section]
\newtheorem{lemma}[theorem]{Lemma}
\newtheorem{proposition}[theorem]{Proposition}

\newtheorem{definition}[theorem]{Definition}

\newtheorem{problem}[theorem]{Problem}

\renewcommand{\arraystretch}{1.2}     % so lines in tables are not crowded
\renewcommand{\geq}{\geqslant}
\renewcommand{\leq}{\leqslant}
\renewcommand{\ge}{\geqslant}

\begin{document}

\renewcommand\atop[2]{\genfrac{}{}{0pt}{}{#1}{#2}}
\newcommand{\Ann}{\textup{Ann}}
\newcommand{\Aut}{\textup{Aut}}
\newcommand{\C}{\textup{C}}
\newcommand{\Char}{\textup{char}}
\newcommand{\eps}{\varepsilon}
\newcommand{\F}{\mathbb{F}}
\newcommand{\GL}{\textup{GL}}
\newcommand{\im}{\textup{im\,}}
\newcommand{\JCF}{\textup{JCF}}
\newcommand{\M}{\textup{M}}
\newcommand{\m}{$\kern1pt--\kern-1pt$}
\newcommand{\N}{\mathcal{N}}
\newcommand{\OO}{\textup{O}}
\newcommand{\partn}{\vdash} % partition symbol
\newcommand{\Part}{\mathcal P}
\newcommand\qbinom[2]{\genfrac{[}{]}{0pt}{}{#1}{#2}}
\newcommand{\Q}{\mathbb{Q}}
\newcommand{\SL}{\textup{SL}}
\newcommand{\SO}{\textup{SO}}
\newcommand{\Sp}{\textup{Sp}}
\newcommand{\T}{\textup{T}}
\newcommand{\U}{\textup{U}}
\newcommand{\w}{\omega}
\newcommand{\Z}{\mathbb{Z}}

\hyphenation{}

\title[Decomposing modular tensor products]{Decomposing modular tensor products, \\and periodicity of `Jordan partitions'}
\author{S.\,P. Glasby, Cheryl E. Praeger, and Binzhou Xia}

\address[Glasby]{
Centre for Mathematics of Symmetry and Computation\\
University of Western Australia\\
35 Stirling Highway\\
Crawley 6009, Australia.\newline Email: {\tt GlasbyS@gmail.com; WWW: \href{http://www.maths.uwa.edu.au/~glasby/}{http://www.maths.uwa.edu.au/$\sim$glasby/}}\newline
Also affiliated with Department of Mathematics, University of Canberra, Australia. }.
\address[Praeger]
{Centre for Mathematics of Symmetry and Computation\\
University of Western Australia\\
35 Stirling Highway\\
Crawley 6009, Australia.\newline Email: {\tt Cheryl.Praeger@uwa.edu.au; WWW: \href{http://www.maths.uwa.edu.au/~praeger}{http://www.maths.uwa.edu.au/$\sim$praeger} }\newline Also affiliated with King Abdulaziz University, Jeddah,
 Saudi Arabia.}
\address[Xia]{
Beijing International Center for mathematical Research\\
Peking University\\
Beijing 100871, People's Republic of China.\newline Email: {\tt binzhouxia@pku.edu.cn}
}

\date{\today}

\begin{abstract}
Let $J_r$ denote an $r\times r$ matrix with minimal and characteristic
polynomials $(t-1)^r$. Suppose $r\leq s$. It is not hard to show
that the Jordan canonical form of $J_r\otimes J_s$ is similar
to $J_{\lambda_1}\oplus\cdots\oplus J_{\lambda_r}$ where
$\lambda_1\geq\cdots\geq\lambda_r>0$ and $\sum_{i=1}^r\lambda_i=rs$.
The partition $\lambda(r,s,p):=(\lambda_1,\dots,\lambda_r)$ of $rs$, which
depends only on $r,s$ and the characteristic $p:=\Char(F)$,
has many applications including to the study of algebraic groups.
We prove new periodicity and duality results for $\lambda(r,s,p)$ that
depend on the smallest $p$-power exceeding~$r$.
This generalizes results of J.\,A.~Green,  B. Srinivasan, and others which
depend on the smallest $p$-power exceeding the (potentially large) integer~$s$.
It also implies that for fixed $r$ we can construct a finite table allowing the
computation of $\lambda(r,s,p)$ for all $s$ and $p$, with $s\geq r$ and
$p$ prime.
%This generalizes work of K-i. Iima and R. Iwamatsu.
%Let $J_r$ denote an $r\times r$ matrix with minimal and characteristic
%polynomials $(t-1)^r$. If $r\leq s$, then the Jordan canonical form of
%$J_r\otimes J_s$ is similar to $J_{\lambda_1}\oplus\cdots\oplus J_{\lambda_r}$
%where $(\lambda_1,\dots,\lambda_r)$ is a partition of $rs$. This partition,
%denoted $\lambda(r,s,p)$ where $p:=\Char(F)$, has many applications including
%to the study of algebraic groups. We show for fixed $r$ that a finite
%computation will find $\lambda(r,s,p)$ for all $s\geq r$, and all primes~$p$.
%(For a precise statement see Theorem~\ref{T:finite}.)
%We also generalize the main result of B. Srinivasan, The modular
%representation ring of a cyclic $p$-group, {\it Proc. London Math. Soc.}
%{\bf 14} (1964) 677--688, and we prove new periodicity and duality results
%for $\lambda(r,s,p)$ that depend on the smallest $p$-power exceeding~$r$.
\end{abstract}

\maketitle
\centerline{\noindent AMS Subject Classification (2010): 15A69, 15A21, 13C05}

\section{Introduction}\label{S:Intro}

Consider a matrix whose minimal and characteristic polynomials equal
$(t-1)^r$. To be explicit, take the $r\times r$ matrix~$J_r$  with
1s in positions $(i,i)$ for $1\leq i\leq r$, and $(i,i+1)$ for $1\leq i<r$,
and zeros elsewhere. Suppose $1\leq r\leq s$. Then the Jordan canonical form of
$J_r\otimes J_s$ is a direct sum $J_{\lambda_1}\oplus\cdots\oplus J_{\lambda_r}$,
with precisely $r$ nonempty blocks, see Lemma~\ref{L:baby}(a).
This decomposition depends on the characteristic~$p$ of the
underlying field\footnote{We may assume that $F=\F_p$ or $\Q$ as the Jordan form of $J_r\otimes J_s$ is invariant under field extensions.}~$F$, and it
determines a partition $\lambda(r,s,p)=(\lambda_1,\dots,\lambda_r)$
of $rs$ since $J_r\otimes J_s$ is an $rs\times rs$ matrix.
We will assume that $\lambda_1\geq\cdots\geq\lambda_r>0$.
The determination of this
`Jordan partition'\footnote{This phrase was used by Dmitri Panyushev in the review MR2728146, but it is not used commonly.} has applications to many significant problems.
The representation theory of algebraic groups is
governed by the behaviour of the unipotent elements, and indeed
properties of $\lambda(r,s,p)$ are particularly useful (when $p>0$) for the
study of exceptional algebraic groups, see~\cite{LS,rL}. More generally,
Lindsey~\cite[Theorem\;1]{jL} gives a useful (though somewhat technical)
lower bound
on the degree of the minimal faithful representation in characteristic~$p$
for certain groups with a prescribed Sylow $p$-subgroup structure.
Lindsey's result, in turn, may be applied to the study of primitive permutation
groups of $p$-power degree, see~\cite{P}.

The most direct application, and the oldest, is to the study of
modular representations of finite cyclic $p$-groups. Given two indecomposable
modules $V_r$ and $V_s$ of a cyclic group $G$ of order $p^n$, the
module $V_r\otimes V_s$ is, by the Krull-Schmidt theorem, a sum of
indecomposable modules $V_{\lambda_1}\oplus\cdots\oplus V_{\lambda_r}$. Thus
when $p>0$, the partition $\lambda(r,s,p)$ arises naturally in this context too.
The connection with matrices is straightforward:
$G=\langle g\rangle$ has precisely~$p^n$ pairwise nonisomorphic
indecomposable modules $V_1,\dots,V_{p^n}$ which correspond to
the matrix representations $G\to\GL(r,\F_p)\colon g\mapsto J_r$
where $1\leq r\leq p^n$.

\begin{definition}\label{D:}
The following terminology will be used as convenient abbreviations.
\begin{itemize}
\item[\textup{(a)}] For integers $r,s$ with $1\leq r\leq s$, the
  \emph{\bf standard partition}
  $\lambda=(\lambda_1,\dots,\lambda_r)$ of $rs$ is the partition with
  $\lambda_i=r+s-2i+1$ for $1\leq i\leq r$, i.e. $(s+r-1,\dots,s-r+1)$.
\item[\textup{(b)}] Call $\lambda=(\lambda_1,\dots,\lambda_r)$ the
  \emph{\bf ($r$-)uniform partition} of $rs$ if $\lambda_i=s$, for
  $1\leq i\leq r$.
\item[\textup{(c)}] The vector
  $\eps(r,s,p)=(\eps_1,\dots,\eps_r)$ with $\eps_i=\lambda_i-s$, which
  measures the deviation of $\lambda(r,s,p)=(\lambda_1,\dots,\lambda_r)$
  from the uniform vector, is called the \emph{\bf deviation vector}.
\item[\textup{(d)}] The \emph{\bf negative reverse} of $(\eps_1,\eps_2,\dots,\eps_r)$
  is $\overline{(\eps_1,\eps_2,\dots,\eps_r)}:=(-\eps_r,\dots,-\eps_2,-\eps_1)$.
\item[\textup{(e)}] The \emph{\bf $k$-multiple} of
  $(\lambda_1,\dots,\lambda_r)$ is the vector $(k\lambda_1,\dots,k\lambda_1,\dots,k\lambda_r,\dots,k\lambda_r)$ of length $kr$ where the size,
and multiplicity, of each part is multiplied by~$k$.
%\item[\textup{(e)}] The \emph{\bf $k$-multiple} of a vector
%  $(\lambda_1,\dots,\lambda_r)$ with $r$ coordinates is the vector
%  with $kr$ coordinates and each $\lambda_i$ repeated $k$ times, \emph{viz.}
%  $(\lambda_1,\dots,\lambda_1,\dots,\lambda_r,\dots,\lambda_r)$.
\end{itemize}
\end{definition}

In characteristic zero, the partition $\lambda(r,s,0)$ was shown to be
the standard partition independently by Aitken (1934), Roth (1934),
and Littlewood (1936); for more background and references see~\cite[p.\,416]{N2}.
The change-of-basis matrix exhibiting
the Jordan canonical form of $J_r\otimes J_s$ may be chosen to have rational
entries, and so in `large' prime characteristic ($p$ not dividing denominators
of the matrix entries), it follows that $\lambda(r,s,p)$ is also the
standard partition. Srinivasan proved that $\lambda(r,s,p)$ is
standard for $p\geq r+s-1$. Our first result generalizes the main
results of both~\cite{S} and~\cite{B}.

\begin{theorem}\label{T:std}
If $r\leq s$, and
$s\not\equiv 0,\pm1,\pm2,\dots,\pm(r-2)\pmod p$, then $\lambda(r,s,p)$ is the
standard partition; i.e. its $i$\emph{th} part is
$\lambda_i=r+s+1-2i$ for $1\leq i\leq r$.
\end{theorem}

Throughout this paper we have
$p>0$ and $r\leq s$. It is useful (psychologically) to think
of~$r$ as `fixed and small', and $s$ and $p$ as `variable', and $s$ as `large'.
The seminal paper~\cite{G} by J.\,A.~Green led to a series of different
algorithms~\cite{R,S,M1,M2,R1,R2,H,II} for decomposing
$V_r\otimes V_s$ as $V_{\lambda_1}\oplus\cdots\oplus V_{\lambda_r}$ where
$1\leq r\leq s\leq p^n$, and $p^n$ is the smallest $p$-power exceeding $s$.
One class of algorithms~\cite{G,S,M1,R1,R2} involves recursive computations
in the
modular representation ring (or Green ring) of the cyclic group $C_{p^n}$.
It is known that $V_{p^i}$ for $0\leq i\leq n$ are generators for the
Green ring; however, relations in these generators are rather mysterious.
Another class of algorithms is related to $p$-adic expansions, and has a more
number-theoretic flavour. The algorithms may depend on $p$-adic
expansions of $r$ and $s$ as in~\cite{M1}, or on the values of certain
determinants modulo~$p$ as in~\cite{S,II}, or on so called
$p$-ranks~\cite{N1,H}.
Ideally the algorithms can construct a basis relative to which
$J_r\otimes J_s$ is in Jordan canonical form, see~\cite{B,N2}.
No complexity analysis presently exists to compare the time or space
requirements of these algorithms.

In contrast to the Green ring results which assume $1\leq r\leq s\leq p^n$, we
assume only that $r\leq\min\{s,p^m\}$. For maximum impact
take $p^m$ to be the \emph{smallest} $p$-power exceeding~$r$. We first noticed
a criterion for $\lambda(r,s,p)$ to be the uniform partition.

\begin{proposition}\label{P:}
If $\Char(F)=p>0$ and $r\leq\min\{p^m,s\}$, then $\lambda(r,s,p)$ is the
uniform partition, or equivalently
$\eps(r,s,p)=(0,\dots,0)$, if and only if $s\equiv 0\pmod{p^m}$.
\end{proposition}

Proposition~\ref{P:}, which is more general than~\cite[Lemma~2.1]{R1},
for example, suggested to us to study the deviation
$\eps(r,s,p)=(\lambda_1,\dots,\lambda_r)-(s,\dots,s)$
from the uniform partition $(s,\dots,s)$. We found several properties
of $\eps(r,s,p)$ that depend only on the congruence of $s$ modulo~$p^m$,
where $m=\lceil\log_p(r)\rceil$, see
Theorem~\ref{T:DP}. These properties generalize periodicity and duality
results for Green rings, \emph{c.f.}~\cite[Eq.\;G-1]{R1}.

\begin{theorem}\label{T:DP}
Suppose $r\leq \min\{s,s',p^m\}$ where $p=\Char(F)$.
\begin{itemize}
  \item[(a)] \textup{[Periodicity]} If $s\equiv s'\pmod{p^m}$, then
  $\eps(r,s,p)=\eps(r,s',p)$.
  \item[(b)] \textup{[Duality]} If $s'\equiv -s\pmod{p^m}$, then $\eps(r,s',p)$
  is the negative reverse of $\eps(r,s,p)$.
\end{itemize}
\end{theorem}

Barry~\cite{B2} has already used Theorem~\ref{T:DP} to classify all triples
$(r,s,p)$ for which $\lambda(r,s,p)$ is standard, and Barry's result
already has an application, namely determining whether so called
`Jordan permutations' are involutions or are trivial, see~\cite{GPX2}.

Section~\ref{S:Basic} introduces notation and terminology whilst establishing an
important result that is computationally advantageous:
computing with the nilpotent matrix $J_r\otimes I_s+ I_r\otimes J_s$
rather than the natural
unipotent matrix $J_r\otimes J_s$, {\it c.f.}~\cite{T1,T2}.
It is convenient to view the tensor product
of vector spaces as a polynomial algebra; then submodules correspond to ideals.
In Section~\ref{S:std} we prove the sufficient condition Theorem~\ref{T:std}
for $\lambda(r,s,p)$ to be standard, and also prove:

\begin{theorem}\label{T:p}
If $r\leq s$ and $k\geq0$, then $\lambda(p^kr,p^ks,p)$ is the
$p^k$-multiple of $\lambda(r,s,p)$.
\end{theorem}

\noindent
Renaud~\cite[Lemma~2.2]{R1} proved this result (using the language of Green
rings)
under the additional assumption that $r\leq s\leq p$ which we do not need.
Our $s$ can be large.

In Section~\ref{S:PD} we prove Proposition~\ref{P:} and Theorem~\ref{T:DP},
and find $\eps(r,s,p)$ explicitly when $s\equiv\pm1$ or $\pm2\pmod{p^m}$, see
Propositions~\ref{P:s1} and~\ref{P:s2}. In Section~\ref{S:small} we prove the following theorem.

\begin{theorem}\label{T:nr}
For a fixed integer $r\geq1$, there are at most $2^{r-1}$
different deviation vectors $\eps(r,s,p)$ as both $s$, where $s\geq r$,
and the prime $p$ vary.
\end{theorem}

It was shown in~\cite[Theorem~2.1.5]{II} that a finite computation is
required to determine $\lambda(r,s,p)$ when both $r$ and $s$ are fixed,
and~$p$~varies. Table~\ref{T:eps} lists the possible deviation vectors
$\eps(r,s,p)$ for $r\leq5$. In Section~\ref{S:small} we prove the
generalization (Theorem~\ref{T:finite}) which allows us to create
Table~\ref{T:eps}.
In the proof of Theorem~\ref{T:finite} we show that, for a given $r$,
the number of values of $s$ and of $p$ we need to consider are each
bounded above in terms of~$r$.

How is Table~\ref{T:eps} used?
This table lists the values of $\eps(r,s,p)$ for $r\leq\min\{5,s\}$.
It explicitly lists the `small' primes $p<2r-3$; these may have $m>1$. The
infinitely many `large' primes $p'\geq 2r-3$ all have $m=1$.
For the small primes $p$ it suffices, by duality, to list the ${s\pmod{p^m}}$
for which $0\leq {s\pmod{p^m}}\leq p^m/2$. For the large primes $p'$
it suffices, by Theorem~\ref{T:std} and duality, to list the ${s\pmod{p'}}$
for which $0\leq {s\pmod{p'}}\leq r-2$. The values of $s\pmod{p'}$
satisfying $r-1\leq s\pmod {p'}\leq p'-(r-1)$, have
$\lambda(r,s,p')$ standard. Thus
$\eps(r,s,p')=(r-1,r-3,\dots,-(r-3),-(r-1))$ for these $p'-(2r-3)$ choices of
${s\pmod{p'}}$. We also list this `standard vector' $\eps(r,s,p')$
for ${s\pmod{p'}}=r-1$. Note that for $(r,p')=(3,3), (4,5), (5,7)$, no value
of $s$ gives rise to this standard vector as $p'-(2r-3)=0$.
When computing $\eps(4,17,3)$ we have $m=2$ as
$3<4\leq 3^2$. The first equality below is by periodicity, the
second is by duality, the third is by Table~\ref{T:eps}, and the fourth is by
the definition (see~\ref{D:}(e)) of `negative reverse':
\[
  \eps(4,17,3)=\eps(4,8,3)=\overline{\eps(4,10,3)}=\overline{(3,-1,-1,-1)}=(1,1,1,-3).
\]
Also $\eps(5,s,11)=(4,2,0,-2,-4)$ for $s=15,16,17,18$, and $p'-(2r-3)=11-7=4$.

\begin{theorem}\label{T:finite}
For fixed $r$, a finite computation suffices to compute the values of
$\eps(r,s,p)$ for all $s$ with $s\geq r$, and all primes $p$.
\end{theorem}

\begin{table}[!ht]
\caption{Values of $\eps(r,s,p)$ with $r\leq\min\{5,s\}$ and $m=\lceil\log_p(r)\rceil$. See the paragraph beginning `How is Table~\ref{T:eps} used?';
below $p'$ is a prime~$\ge 2r-3$.}\label{T:eps}
%for explication of use.}\label{T:eps}
\begin{tabular}{ |c|l|l|l|l|l| } \hline
$s\kern-2pt\mod{p^m}$&0&1&2&3&4\\ \hline
$\eps(1,s,p')$ & (0) & & &&\\ \hline
$\eps(2,s,p')$ & (0,0) & $(1,\m1)$& &&\\ \hline
$\eps(3,s,2)$ & (0,0,0) & $(2,\m1,\m1)$&$(2,0,\m2)$&&\\ %\hline
$\eps(3,s,p')$ & (0,0,0) & $(2,\m1,\m1)$&$(2,0,\m2)$&&\\ \hline
$\eps(4,s,2)$ & (0,0,0,0) & $(3,\m1,\m1,\m1)$&$(2,2,\m2,\m2)$&&\\ %\hline
$\eps(4,s,3)$ & (0,0,0,0) & $(3,\m1,\m1,\m1)$&$(3,1,\m2,\m2)$&$(3,0,0,\m3)$&$(3,1,\m1,\m3)$\\ %\hline
$\eps(4,s,p')$ & (0,0,0,0) & $(3,\m1,\m1,\m1)$&$(3,1,\m2,\m2)$&$(3,1,\m1,\m3)$&\\ \hline
$\eps(5,s,2)$ & (0,0,0,0,0) & $(4,\m1,\m1,\m1,\m1)$&$(4,2,\m2,\m2,\m2)$&$(4,1,1,\m3,\m3)$&$(4,0,0,0,\m4)$\\ %\hline
$\eps(5,s,3)$ & (0,0,0,0,0) & $(4,\m1,\m1,\m1,\m1)$&$(4,2,\m2,\m2,\m2)$&$(3,3,\kern-1pt0,\kern-1pt\m3,\kern-1pt\m3)$&$(4,\kern-1pt2,\kern-1pt0,\kern-1pt\m2,\kern-1pt\m4)$\\ %\hline
$\eps(5,s,5)$ & (0,0,0,0,0) & $(4,\m1,\m1,\m1,\m1)$&$(3,3,\m2,\m2,\m2)$& &\\ %\hline
%$\eps(5,s,7)$ & (0,0,0,0,0) & $(4,\m1,\m1,\m1,\m1)$&$(4,2,\m2,\m2,\m2)$&$(4,2,0,\m3,\m3)$&\\ %\hline
$\eps(5,s,p')$ & (0,0,0,0,0) & $(4,\m1,\m1,\m1,\m1)$&$(4,2,\m2,\m2,\m2)$&$(4,2,0,\m3,\m3)$&$(4,\kern-1pt2,\kern-1pt0,\kern-1pt\m2,\kern-1pt\m4)$\\ \hline
\end{tabular}
\end{table}

\section{Notation and basic results}\label{S:Basic}

This section introduces notation and establishes facts needed for
proofs in subsequent sections. Parts (a--c) of Lemma~\ref{L:baby}
have been proved
before in~\cite[Lemma\;2.1]{R}, \cite[p.\;678]{S}, and \cite[Theorem~2]{M2},
but because we want to build on their proofs,
it is desirable to give new proofs using our polynomial notation.
Our alternative proofs, which are based on Lemma~\ref{L:coeff} and the preamble
to Lemma~\ref{L:baby}, are much shorter than the original proofs.

Fix a positive integer $r$ and a field $F$, and consider the quotient polynomial
ring $B:=F[X]/(X^r)$. Set $x:=X+(X^r)$. Then $x^r=0$
and $1,x,\dots,x^{r-1}$~is a basis for~$B$. As usual, right multiplication gives
rise to a monomorphism $\mu\colon B\to\textup{End}_F(B)$ where for
$b\in B$ the $F$-linear map $\mu_b:=\mu(b)$ satisfies $\mu_b(a)=ab$.
Denote the matrices of $\mu_x$, $\mu_{1+x}$, and $\mu_{\alpha+x}$
relative to the basis $1,x,\dots,x^{r-1}$~by
\[
  N_r=\begin{pmatrix}0&1&&\\&\ddots&\ddots&\\&&0&1\\&&&0\end{pmatrix},\quad
  J_r=\begin{pmatrix}1&1&&\\&\ddots&\ddots&\\&&1&1\\&&&1\end{pmatrix},\quad
  J_r(\alpha)=\begin{pmatrix}\alpha&1&&\\&\ddots&\ddots&\\&&\alpha&1\\&&&\alpha\end{pmatrix},
\]
respectively. We say $N_r$ is {\it nilpotent}, $J_r$ is {\it unipotent}, and $J_r(\alpha)$ is {\it $\alpha$-potent}, i.e. $J_r(\alpha)-\alpha I$ is nilpotent. Moreover, $J_r(\alpha)$ is an $r\times r$ Jordan block,
and its minimal polynomial is $(t-\alpha)^r$. The matrix
$J_r(\alpha)\otimes J_s(\beta)$ is an upper-triangular $\alpha\beta$-potent
matrix, and so its Jordan canonical form is
$\JCF(J_r(\alpha)\otimes J_s(\beta))=\bigoplus_{t\geq1} c_{r,s,t}J_t(\alpha\beta)$
where $c_{r,s,t}\in\mathbb{N}$ denotes the multiplicity of the $t\times t$ Jordan
block $J_t(\alpha\beta)$. The Jordan partition
$\lambda:=\langle t^{c_{r,s,t}}\rangle$ of~$rs$, with part size~$t$ occurring with
multiplicity $c_{r,s,t}$, has been studied by many authors. In the
nilpotent case (when $\alpha\beta=0$) the multiplicities $c_{r,s,t}$
are easily described (see for example~\cite[2.1.2]{II}). Furthermore,
the change of basis matrix is known, and is field independent.
The invertible case (when $\alpha\beta\neq0$) reduces to the unipotent case
because
\[
  \JCF(J_r(\alpha)\otimes J_s(\beta))=\bigoplus_{t\geq1} c_{r,s,t}J_t(\alpha\beta)
  \quad\textup{if and only if}\qquad
  \JCF(J_r\otimes J_s)=\bigoplus_{t\geq1} c_{r,s,t}J_t.
\]
Thus the same Jordan partitions arise, and the change of basis matrices
are easily related.

The quotient polynomial algebra $A:=F[X,Y]/(X^r,Y^s)$ is an $rs$-dimensional
$F$-vector space with basis $x^iy^j$, $0\leq i<r$, $0\leq j<s$, where
$x:=X+(X^r,Y^s)$ and $y:=Y+(X^r,Y^s)$ satisfy $x^r=y^s=0$.
Given $a\in A$ denote by $F[a]$ the subalgebra
$\{f(a)\mid f(t)\in F[t]\}$ of $A$. Then $A$, viewed as a module over the
ring $F[a]$, is a direct sum $A=a_1F[a]\oplus\cdots\oplus a_nF[a]$ of
cyclic $F[a]$-submodules. To avoid ambiguity, we regard $A$
as an $F[a]$-module rather than more conventionally as an $F[t]$-module,
see~\cite{HH}. Indeed, $F[a]$ is a quotient of the principal ideal domain
$F[t]$. We are interested in the dimensions of the
cyclic submodules when $a=(1+x)(1+y)$. Since $F[a]=F[a-1]$ it is convenient to
replace the invertible unipotent element $(1+x)(1+y)$ with the nilpotent
element $(1+x)(1+y)-1=x+y+xy$. We show that $x+y+xy$
and $x+y$ induce {\it similar} linear transformations on $A$.
The action of $x+y$ on the basis $x^iy^j$ is simple:
$x^iy^j(x+y)=x^{i+1}y^j+x^iy^{j+1}$. We seek another basis $f_{i,j}$ for $A$
such that $f_{i,j}(x+y+xy)=f_{i+1,j}+f_{i,j+1}$. An easy calculation shows that
$f_{i,j}=x^i(1+y)^iy^j$, $0\leq i<r$, $0\leq j<s$, is the desired basis.
This is a major point in~\cite{T1,T2}.
View $A$ as a module over
\begin{equation}\label{E:pol}
  F[(1+x)(1+y)]=F[x+y+xy], \quad\textup{or over}\quad F[x+y].
\end{equation}
(Incidentally, when decomposing
$J_r\otimes J_s\otimes J_t$ it is similarly useful to consider a module
over $F[x+y+z]$ instead of $F[(1+x)(1+y)(1+z)]$ by using $f_{i,j,k}=x^i(1+y)^i(1+z)^iy^j(1+z)^jz^k$.)

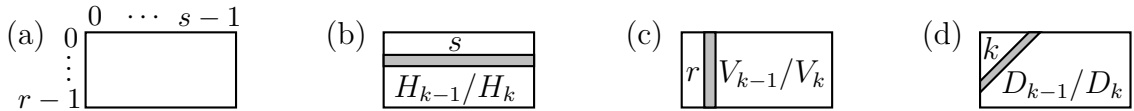
\begin{figure}[!ht]
  \caption{(a) $r\times s$ basis elements $x^iy^j$;
    and (b)--(d) dimensions of sections.}
  \vskip5mm
  \label{F:dim}
  \begin{center}
  \psset{xunit=5mm}\psset{yunit=5mm}
  \begin{pspicture}(-2.6,0)(5,2)
    \thicklines
    \pspolygon(0,0)(4,0)(4,2)(0,2)
    \uput[180](0,2){\rput(-0.04,-0.15){{\small $0$}}}
    \uput[180](0,1){\rput(-0.1,0.25){{\small $\vdots$}}}
    \uput[180](0,0){\rput(-0.6,0.2){{\small $r-1$}}}
    \uput[90](0,2){\rput(0.25,0.1){{\small $0$}}}
    \uput[90](1.5,2){\rput(0.1,0.05){{\small $\cdots$}}}
    \uput[90](3,2){\rput(0.3,0.05){{\small $s-1$}}}
    \rput(-1.6,1.9){(a)}
  \end{pspicture}
  \qquad
    \begin{pspicture}(-1,0)(5,2)
    \thicklines
    \pspolygon(0,0)(4,0)(4,2)(0,2)
    \pspolygon[fillstyle=solid,fillcolor=lightgray](0,1.1)(0,1.4)(4,1.4)(4,1.1)
    \rput(-1,1.9){(b)}    \rput(1.9,0.5){$H_{k-1}/H_k$}\rput(1.9,1.7){$s$}
  \end{pspicture}
  \qquad
    \begin{pspicture}(-1,0)(5,2)
    \thicklines
    \pspolygon(0,0)(4,0)(4,2)(0,2)
    \pspolygon[fillstyle=solid,fillcolor=lightgray](0.6,0)(0.9,0)(0.9,2)(0.6,2)
    \rput(-1,1.9){(c)}    \rput(2.4,0.9){$V_{k-1}/V_k$}\rput(0.3,0.9){$r$}
  \end{pspicture}
  \qquad
    \begin{pspicture}(-1,0)(5,2)
    \thicklines
    \pspolygon(0,0)(4,0)(4,2)(0,2)
    \pspolygon[fillstyle=solid,fillcolor=lightgray](0,0.4)(1.6,2)(1.3,2)(0,0.7)
    \rput(-1,1.9){(d)}    \rput(2.2,0.6){$D_{k-1}/D_k$}\rput(0.3,1.6){$k$}
  \end{pspicture}
  \end{center}
\end{figure}

View the basis elements $x^iy^j$, $0\leq i<r$, $0\leq j<s$,
of $A$ as placed on an $r\times s$ rectangle, see Figure~\ref{F:dim}(a).
Consider the
following `horizontal', `vertical', and `diagonal' ideals of $A$:
\begin{align*}
  H_k&=\langle x^iy^j\mid i\geq k\rangle, &&A=H_0>H_1>\cdots>H_r=0,\\
  V_k&=\langle x^iy^j\mid j\geq k\rangle, &&A=V_0>V_1>\cdots>V_s=0,\\
  D_k&=\langle x^iy^j\mid i+j\geq k\rangle, &&A=D_0>D_1>\cdots>D_{r+s-1}=0.
\end{align*}

The matrix of $\mu_{(1+x)(1+y)}$ is $J_r\otimes J_s$ relative to the basis
$x^iy^j$, $0\leq i<r$, $0\leq j<s$, ordered lexicographically by $i$, then $j$.
The action of $\mu_{(1+x)(1+y)}$ on the above ideals and their quotients
(relative to this monomial basis) is given in Table~\ref{T:Action}.

\begin{table}[!ht]
\renewcommand{\arraystretch}{1.5}     % so lines in tables are not crowded
\caption{The action of $\mu_{(1+x)(1+y)}$ on submodules
 and quotient modules of $A$.}
\label{T:Action}\kern-2mm\begin{center}
\begin{tabular}{|c|c|c|c|c|c|c|} \hline
$H_k$&$A/H_k$&$V_k$&$A/V_k$&$H_{k-1}/H_k$&$V_{k-1}/V_k$&$D_{k-1}/D_k$\\ \hline
$J_{r-k}\otimes J_s$&$J_k\otimes J_s$&$J_r\otimes J_{s-k}$&$J_r\otimes J_k$&$J_s$&$J_r$&$I$\\ \hline
\end{tabular}
\end{center}
\end{table}

Since $H_{k-1}/H_k$ is a cyclic $F[x+y]$-module (generated by $x^{k-1}+H_k$),
it follows that $A$ is a sum of at most $r$ cyclic $F[x+y]$-submodules.
Indeed, $A=\sum_{i=0}^{r-1} x^iF[x+y]$.
Conversely, $A$ is a sum of
at least~$r$ cyclic $F[x+y]$-submodules because $D_{r-1}/D_r$ is
an $r$-dimensional vector space,
and $D_{r-1}(x+y)\subseteq D_r$.
Thus $A$ is a sum of precisely $r$ nonzero cyclic $F[x+y]$-submodules and
$\lambda(r,s,p)$ has precisely~$r$ nonzero parts,
see Lemma~\ref{L:baby}(a) below.
The dimensions of the sections $H_{k-1}/H_k$, $V_{k-1}/V_k$ and $D_{k-1}/D_k$
can be seen from Figure~\ref{F:dim}(b)--(d) to be:
\begin{align*}
  &\dim(H_{k-1}/H_k)=s                &&\textup{for $1\leq k\leq r$,}\\
  &\dim(V_{k-1}/V_k)=r                &&\textup{for $1\leq k\leq s$,}\\
  &\dim(D_{k-1}/D_k)=\min\{k,r,r+s-k\}&&\textup{for $1\leq k<r+s$.}
\end{align*}
Define the {\it annihilator} of an element $a\in A$ to be
$\Ann(a)=\{b\in A\mid ab=0\}$. Then $\Ann(x+y)$ is an $r$-dimensional
$F[x+y]$-submodule of $A$ because $\lambda(r,s,p)$ has
precisely~$r$ nonzero parts, as proved above; see also Lemmas~\ref{L:coeff}(b)
and~\ref{L:baby}(a).

Henceforth view $A$ as an $F[x+y]$-module. The dimension of the
cyclic submodule $aF[x+y]$ is the smallest natural number $n=n(a)$
satisfying $a(x+y)^n=0$.
Clearly $(x+y)^i\neq0$ for $0\leq i<s$. Simplifying
$(x+y)^n=\sum_{i+j=n} \binom{n}{i}x^iy^{j}$ using $x^r=y^s=0$ gives
\begin{equation}\label{E:bin}
  (x+y)^n=\sum_{i=n-s+1}^{r-1} \binom{n}{i}x^iy^{n-i}
  \qquad\qquad\textup{for $s\leq n< r+s-1$.}
\end{equation}
Eq.~\eqref{E:bin} implies that $(x+y)^n=0$ for $n\geq r+s-1$ because an empty
sum is zero. This shows that $s\leq\lambda_1\leq r+s-1$.
Homogeneous polynomials will play an important role.

\begin{lemma}\label{L:coeff}
Suppose $\Char(F)=p$ and $w=\sum_{i=n-s+1}^{r-1}\alpha_ix^iy^{n-i}$ where
$s-1\leq n\leq r+s-2$.
\begin{itemize}
\item[(a)] $w(x+y)=0$ holds if and only if $w$ is a scalar multiple of $\sum_{i=n-s+1}^{r-1}(-1)^ix^iy^{n-i}$.
\item[(b)] If $w_i:=\sum_{j=0}^i(-1)^jx^{r-1-j}y^{s-1-i+j}$, then $w_0,\dots,w_{r-1}$
is an $F$-basis for $\Ann(x+y)$.
\item[(c)] If $\dim(F[x+y])=n$, then $(x+y)^n=0$ and $(x+y)^{n-1}\neq0$.
Moreover,
\[
  (x+y)^{n-1}=(-1)^{r-1}\binom{n-1}{r-1}\sum_{i=n-s}^{r-1} (-x)^iy^{n-1-i}
  \quad\textup{and $\binom{n-1}{r-1}\neq0$ in $F$.}
\]
\end{itemize}
\end{lemma}

\begin{proof}
(a) Expanding $w(x+y)$ and using $x^r=y^s=0$ gives
\begin{align*}
  &\alpha_{n-s+1}x^{n-s+2}y^{s-1}+\alpha_{n-s+2}x^{n-s+3}y^{s-2}+\cdots+
           \alpha_{r-3}x^{r-2}y^{n-r+3}+\alpha_{r-2}x^{r-1}y^{n-r+2}+\\
  &\alpha_{n-s+2}x^{n-s+2}y^{s-1}+\alpha_{n-s+3}x^{n-s+3}y^{s-2}+\cdots+
           \alpha_{r-2}x^{r-2}y^{n-r+3}+\alpha_{r-1}x^{r-1}y^{n-r+2}.
\end{align*}
Thus $w(x+y)=0$ holds if and only if $\alpha_i=-\alpha_{i+1}$
for $n-s+1\leq i\leq r-2$, as desired.

(b) As $w_i$ is homogeneous of degree $r+s-2-i$, it follows that
$w_0,\dots,w_{r-1}$ are $F$-linearly independent. The proof of part~(a)
shows that $w_i(x+y)=0$ for each~$i$. This proves that $\dim(\Ann(x+y))\geq r$.
To prove $\dim(\Ann(x+y))\leq r$, we relate $v$ and $v':=v(x+y)$. If
\[
  v=\sum_{i=0}^{r-1}\sum_{j=0}^{s-1}v_{i,j}x^iy^j\quad\textup{and}\quad
  v'=\sum_{i=0}^{r-1}\sum_{j=0}^{s-1}v'_{i,j}x^iy^j,
\]
then $v'_{i,n+1-i}=v_{i-1,n-(i+1)}+v_{i,n-i}$. Suppose there exists
$n<r+s-2$ and $i<r$ such that
\[
  v_{i',j'}=0\textup{ whenever $i'+j'<n$,}\quad v_{i',n-i'}=0
  \textup{ for $i'<i$, and}\quad v_{i,n-i}\ne0.
\]
Then $v'_{i,n+1-i}=v_{i,n-i}$, so $v'_{i,n+1-i}x^{i}y^{n+1-i}\ne0$ and $v(x+y)\ne0$.
Hence an element $v\in\Ann(x+y)$ is a sum $v=\sum_{n=s-1}^{r+s-2}v_n$, where
each summand is a homogeneous polynomial
$v_n=\sum_{i=n-s+1}^{r-1}v_{i,n-i}x^iy^{n-i}$ of
degree~$n$. Since $v_n(x+y)$ is homogeneous of degree~$n+1$, the equation
$\sum_{n=s-1}^{r+s-2}v_n(x+y)=0$ implies $v_n(x+y)=0$ for $s-1\leq n\leq r+s-2$.
It now follows from part~(a) that $v$ is a linear combination
of $w_0,w_1,\dots,w_{r-1}$. Thus $w_0,w_1,\dots,w_{r-1}$ is
a basis for $\Ann(x+y)$
as claimed.

(c) Now $\dim(F[x+y])=n$ implies $(x+y)^n=0$. Since $(x+y)^{n-1}$ is
a homogeneous polynomial in $\Ann(x+y)$, part~(a) shows that
$(x+y)^{n-1}=\beta\sum_{i=n-s}^{r-1} (-x)^iy^{n-1-i}$ for some $\beta\in F$.
The binomial theorem shows
that $\beta=(-1)^{r-1}\binom{n-1}{r-1}$. Since $(x+y)^{n-1}\neq0$,
we must have $\beta\neq0$ in~$F$.
\end{proof}

It follows from Lemma~\ref{L:coeff}(c) that $\dim(x^0F[x+y])=\lambda_1$,
and the coefficient of $x^{\lambda_1-s}y^{s-1}$ in $(x+y)^{\lambda_1-1}$ is nonzero.
Now $x^i(x^{\lambda_1-s}y^{s-1})\neq0$ for $0\leq i\leq r+s-1-\lambda_1$. Thus
$x^i(x+y)^{\lambda_1-1}\neq0$ and hence
$\dim(x^iF[x+y])=\lambda_1$. Therefore the $F[x+y]$-submodule
$A'=\bigoplus_{i=0}^{r+s-1-\lambda_1} x^iF[x+y]$ has
dimension $\lambda_1(r+s-\lambda_1)$ as each of the $r+s-\lambda_1$ summands
has dimension~$\lambda_1$.
Any $F[t]$-submodule (or $F[x+y]$-submodule) of maximal dimension
$\lambda_1$ is known to have an $F[t]$-submodule
complement, see~\cite[Ex.\;9.2,\;p.\;149]{HH}. Hence
$A=A'\oplus A''$ for some $F[x+y]$-submodule $A''$. Lemma~\ref{L:coeff}(c)
gives $\dim(x^iF[x+y])<\lambda_1$ for $i\geq r+s-\lambda_1$.
Since $A=\sum_{i=0}^{r-1} x^iF[x+y]$, it follows that
$A''(x+y)^{\lambda_1-1}=0$ and the
multiplicity $k$ of the part~$\lambda_1$ is $k=r+s-\lambda_1$. Substituting
$n=\lambda_1=r+s-k$ into Lemma~\ref{L:coeff}(c) gives an
alternate proof of Lemma~\ref{L:baby}(c) below.
In the next section we prove that Theorem~\ref{T:std}
implies Lemma~\ref{L:baby}(b) when $p>0$, and hence when $p\geq0$.
The ideas established in this section will be needed for
the proof of our main results.

\begin{lemma}\label{L:baby}
Suppose $1\leq r\leq s$ and $F$ is a field of characteristic $p\geq 0$.
\begin{itemize}
\item[(a)] If $p\geq0$, then the partition $\lambda(r,s,p)$ of~$rs$
has precisely $r$ nonzero parts.
\item[(b)] If $p=0$ or $p\geq r+s-1$, then $\lambda_i=r+s-2i+1$ for $1\leq i\leq r$.
\item[(c)] If $p>0$, then $\lambda_1=r+s-k$ where $k\geq1$ is minimal such
  that $p\nmid \binom{r+s-1-k}{r-1}$. Moreover, the part $\lambda_1$ has
  multiplicity $k$ in the partition $\lambda(r,s,p)$.
\end{itemize}
\end{lemma}

\begin{proof}
Parts (a), (b) and (c) have been proved by Ralley, Srinivasan, and McFall,
respectively in~\cite[Lemma\;2.1]{R}, \cite[p.\;678]{S},
and \cite[Theorem~2]{M2}.
\end{proof}

Lemma~\ref{L:baby}(c) suggests an algorithm for computing the largest part
$\lambda_1$ and its multiplicity. It is sometimes difficult to predict the
output of this algorithm, but if $s$ is `not too large' then
$\lambda_1=p^m$ where $m=\lceil\log_p(r)\rceil$. A precise formulation is
given in~\cite{GPX}. More importantly, new symmetries
and applications of these symmetries are described in~\cite{GPX}.

\section{The standard partition}\label{S:std}

Recall that $r\leq s$ and $p>0$ is prime. The hypothesis
${s\not\equiv 0,\pm1,\dots,\pm(r-2)\pmod p}$ in Theorem~\ref{T:std}
implies that $p\geq 2r-3$. As we are thinking of $s$ as `large' compared
to $r$, this is a weaker hypothesis than $p\geq r+s-1$
in~\cite[Corollary\;1]{B}.

\begin{proof}[Proof of Theorem~\textup{\ref{T:std}}]
The strategy of this proof is to show that there exist
elements $v_0,v_1,\dots ,v_{r-1}\in A$ satisfying
\begin{itemize}
  \item[(i)] $\dim(v_iF[x+y])=r+s-2i+1$, and
  \item[(ii)] $v_iF[x+y]\cap \sum_{j\neq i} v_jF[x+y]=0$.
\end{itemize}
The sum $\sum_{i=0}^{r-1} v_iF[x+y]$ is direct by (ii). Since
$\sum_{i=0}^{r-1} (r+s-2i+1) = rs$ holds, it follows by (i) that
$A=\bigoplus_{i=0}^{r-1} v_iF[x+y]$ and $\lambda_i=r+s-2i+1$, as claimed.

Choose an $F$-basis $w_0,w_1,\dots,w_{r-1}$ for $\Ann(x+y)$.
We prove that under the stated hypotheses there exist elements
$v_i\in A$ satisfying $v_i(x+y)^{s+r-2-2i}=w_i$
for each~$i$. This implies that~(i) holds, and (ii) holds because
$w_iF[x+y]$ is the unique minimal $F[x+y]$-submodule of
$v_iF[x+y]$ and $w_iF[x+y]\cap \bigoplus_{j\neq i} w_jF[x+y]=0$ holds.

Set $w_i:=\sum_{j=0}^i(-1)^j x^{r-1-j}y^{s-1-i+j}$. Then $w_0,w_1,\dots,w_{r-1}$
is a basis for $\Ann(x+y)$ by Lemma~\ref{L:coeff}(b).
Suppose that $v_i$ has the form
$v_i=\sum_{k=0}^i\nu_k x^{i-k}y^k$ where the $\nu_k\in F$ are unknowns
which (we will see) can be chosen so that
$v_i(x+y)^{s+r-2-2i}=w_i$. The heart of the proof is that, for each
$i$, we can solve for the $\nu_k$.  Each of $v_i$, $(x+y)^{s+r-2-2i}$
and $w_i$ is homogeneous, and
$\deg(v_i)+\deg((x+y)^{s+r-2-2i})=\deg(w_i)$. The binomial
theorem gives
\begin{align*}
  v_i(x+y)^{s+r-2-2i}&=\left(\sum_{k=0}^i\nu_k x^{i-k}y^k\right)
  \left(\sum_{\ell=0}^{s+r-2-2i}\binom{s+r-2-2i}{\ell}x^\ell y^{s+r-2-2i-\ell}\right)\\
  &=\sum_{k=0}^i\sum_{\ell=0}^{s+r-2-2i}\nu_k\binom{s+r-2-2i}{\ell}x^{i-k+\ell}y^{s+r-2-2i-\ell+k}.
\end{align*}
However, $i-k+\ell\leq r-1$ and $s+r-2-2i-\ell+k\leq s-1$ implies
$r-1-2i+k\leq\ell\leq r-1-i+k$. Setting $j:=(r-1-i+k)-\ell$ gives
the simpler range $i\geq j\geq 0$. Thus
\begin{align*}
  v_i(x+y)^{s+r-2-2i}&=\sum_{k=0}^i\sum_{j=0}^i\,\nu_k\binom{s+r-2-2i}{r-1-i-j+k}x^{r-1-j}y^{s-1-i+j}\\
  &=\sum_{j=0}^i\left[\sum_{k=0}^i\,\nu_k\binom{s+r-2-2i}{s-1-i+j-k}\right]x^{r-1-j}y^{s-1-i+j}
\end{align*}
As we want $v_i(x+y)^{s+r-2-2i}=w_i=\sum_{j=0}^i(-1)^j x^{r-1-j}y^{s-1-i+j}$,
equating coefficients of $x^{r-1-j}y^{s-1-i+j}$ for $0\leq j\leq i$ gives
the linear system
\begin{equation}\label{E:A}
  \begin{pmatrix}\nu_0,\nu_1,\dots,\nu_i\end{pmatrix}A_{i+1}=
  \begin{pmatrix}1,-1,\dots,(-1)^i\end{pmatrix}\quad\textup{where}\quad
  A_{i+1}=\begin{pmatrix}\binom{s+r-2-2i}{s-1-i+j-k}\end{pmatrix}_{0\leq j,k\leq i}.
\end{equation}
The matrix
$A_i=\begin{pmatrix}\binom{s+r-2i}{s-i+j-k}\end{pmatrix}_{0\leq j,k\leq i-1}$
equals the matrix $M(i-1)$ defined on~\cite[p.\;145]{II}.
The determinants $\delta_i:=\det(A_i)$, $1\leq i\leq r-1$, play an
important role. Since $A_r$ is upper triangular, we see $\delta_r=1$.
The following formula for $\delta_i$ is given on~\cite[p.\;145]{II}:
\begin{equation}\label{E:delta1}
  \delta_i=\det(A_i)=\prod_{j=0}^{i-1}
    \frac{\binom{r+s-2i+j}{s-i}}{\binom{s-i+j}{s-i}}
     =\prod_{j=0}^{i-1}\prod_{k=0}^{s-i-1}\frac{r+s-2i+j-k}{s-i+j-k}.
\end{equation}
However, $(r+s-2i+j)-(r-i+k)=s-i+j-k$ is a factor of the numerator and
the denominator for $k=0,1,\dots,s-r-1$. We cancel these factors, and
use the falling factorial notation $n^{\underline{i}}:=n(n-1)\cdots(n-i+1)$
for $i>0$ and $n^{\underline{0}}:=1$. For $0\leq i\leq r-1$ we have
\begin{align}\label{E:delta2}
  \delta_i&=\prod_{j=0}^{i-1}
    \frac{(r+s-2i+j)(r+s-2i+j-1)\cdots(s-i+j+1)}{(r-i+j)(r-i+j-1)\cdots(j+1)}\nonumber\\
  &=\frac{(r+s-i-1)^{\underline{i}}\;(r+s-i-2)^{\underline{i}}\;\cdots\;s^{\underline{i}}}{(r-1)^{\underline{i}}\;(r-2)^{\underline{i}}\;\cdots\;i^{\underline{i}}}\\
  &=\prod_{k=0}^{r-1-i}\frac{(r+s-1-i-k)^{\underline{i}}}{(r-1-k)^{\underline{i}}}.\nonumber
\end{align}
Equation~\eqref{E:delta2} is most helpful when $i$ is close to~$r-1$,
and Equation~\eqref{E:delta1} when $i$ is close to~1.
The following variant of~\eqref{E:delta1} uses the identity
$\binom{r+s-2i+j}{r-i+j}\binom{r-i+j}{j}=\binom{r+s-2i+j}{r-i}\binom{s-i+j}{j}$:
\begin{equation}\label{E:delta3}
  \delta_i=\det(A_i)
    =\prod_{j=0}^{i-1}\frac{\binom{r+s-2i+j}{s-i}}{\binom{s-i+j}{s-i}}
    =\prod_{j=0}^{i-1}\frac{\binom{r+s-2i+j}{r-i+j}}{\binom{s-i+j}{j}}
    =\prod_{j=0}^{i-1}\frac{\binom{r+s-2i+j}{r-i}}{\binom{r-i+j}{j}}.
\end{equation}
\begin{table}[ht!]
\renewcommand{\arraystretch}{1.5}
\caption{Values of $\delta_i$ computed using Eq.~\eqref{E:delta3} for $i$ small, and Eq.~\eqref{E:delta2} for $i$ large.}\label{T:deltavalues}
\begin{tabular}{ |c|c|c|c|c|c|c| } \hline
$\delta_0$&$\delta_1$&$\delta_2$&$\quad\cdots\quad$&$\delta_{r-2}$&$\delta_{r-1}$&$\delta_r$\\ \hline
$1$&$\binom{r+s-2}{r-1}$&$\frac{1}{r-1}\binom{r+s-4}{r-2}\binom{r+s-3}{r-2}$&$\cdots$&$\binom{s+2}{r-1}\binom{s+1}{r-2}$&$\binom{s}{r-1}$&$1$\\ \hline
\end{tabular}
\end{table}

A sufficient condition for
$\lambda(r,s,p)$ to be standard is that
$\delta_{1}\delta_{2}\cdots\delta_{r-1}\not\equiv0\pmod p$.
Consider a lower bound for a typical factor
$(r+s-1-i-k)^{\underline{i}}=\prod_{\ell=0}^{i-1}(r+s-1-i-k-\ell)$
of the numerator of ~\eqref{E:delta2}.
Using $0\leq \ell<i$, $0\leq k\leq r-1-i$, and $1\leq i\leq r-1$ gives
\begin{equation}\label{E:lb}
  r+s-1-i-k-\ell\geq r+s-2i-k\geq s-i+1\geq s-(r-2).
\end{equation}
(Note that a lower bound for a factor of the numerator of ~\eqref{E:delta1} is
too small. The canceling required to deduce~\eqref{E:delta2} from
\eqref{E:delta1} was necessary.) An upper bound for a factor of
the numerator of ~\eqref{E:delta2} can be similarly deduced as follows:
\begin{equation}\label{E:ub}
  r+s-1-i-k-\ell\leq r+s-1-i-k\leq r+s-1-i\leq s+(r-2).
\end{equation}
Equations \eqref{E:lb} and \eqref{E:ub} prove that the factors of
the numerator of $\delta_i$ are bounded between $s-(r-2)$ and $s+(r-2)$.
Hence the assumption $s\not\equiv 0,\pm1,\pm2,\dots,\pm(r-2)\pmod p$
implies that
$\delta_{1}\delta_{2}\cdots\delta_{r-1}\not\equiv0\pmod p$. Therefore
$(\nu_0,\nu_1,\dots,\nu_i)$ in Eq.~\eqref{E:A} can be found, and
$\dim(v_iF[x+y])=r+s-1-2i$ holds by Lemma~\ref{L:coeff}(c).
Thus $\eps(r,s,p)$ is standard.
\end{proof}

We now prove that $\lambda(pr,ps,p)$ is the $p$-multiple
of $\lambda(r,s,p)$, see Definition~\ref{D:}(e).

\begin{proof}[Proof of Theorem~\textup{\ref{T:p}}]
It suffices to prove the result when $k=1$.
Set $F:=\F_p$, and consider the $p^2rs$-dimensional $F$-algebra $\hat{A}$ with
commuting generators $\hat{x}$ and~$\hat{y}$, and relations
$\hat{x}^{pr}=\hat{y}^{ps}=0$.
The $F$-subalgebra $A$ generated by
$x:=\hat{x}^p$ and $y:=\hat{y}^p$ has dimension~$rs$ and satisfies
$x^r=y^s=xy-yx=0$. By the binomial theorem, $(\hat{x}+\hat{y})^p=x+y$.
Suppose that the decomposition of $A$ into cyclic
$F[x+y]$-submodules gives rise to the partition
$\lambda(r,s,p)=(\lambda_1,\dots,\lambda_r)$. To determine
the partition $\lambda(pr,ps,p)$ we consider (using~\eqref{E:pol}) the decomposition of $\hat{A}$
into cyclic $F[\hat{x}+\hat{y}]$-submodules.

Suppose that $A=\bigoplus_{i=1}^r a_iF[x+y]$ where $\dim(a_iF[x+y])=\lambda_i$.
We will prove that
$\hat{A}=\bigoplus_{i=1}^r \bigoplus_{j=0}^{p-1} a_i\hat{x}^jF[x+y]$
where $\dim(a_i\hat{x}^{j}F[\hat{x}+\hat{y}])=p\lambda_i$ for $0\leq j<p$.
Since $a_i(x+y)^{\lambda_i}=0$ and $a_i(x+y)^{\lambda_i-1}\neq 0$, it follows that
$a_i(\hat{x}+\hat{y})^{p\lambda_i}=0$ and $a_i(\hat{x}+\hat{y})^{p(\lambda_i-1)}\neq0$.
However, to show $\dim(a_i\hat{x}^{j}F[\hat{x}+\hat{y}])=p\lambda_i$, we must
prove that $a_i\hat{x}^j(\hat{x}+\hat{y})^{p\lambda_i-1}\neq 0$.

Since $0\neq a_i(x+y)^{\lambda_i-1}\in A$ and $(\hat{x}+\hat{y})^p=x+y$,
there exist scalars
$\alpha_{i,i',j'}\in F$, not all zero, such that
$a_i(\hat{x}+\hat{y})^{p\lambda_i-p}=\sum_{i'=0}^{r-1}\sum_{j'=0}^{s-1}\alpha_{i,i',j'}\hat{x}^{pi'}\hat{y}^{p j'}$.
Since $\binom{p-1}{i''}$ equals $(-1)^{i''}$ in~$\F_p$,
we have
$(\hat{x}+\hat{y})^{p-1}=\sum_{i''=0}^{p-1}(-1)^{i''}\hat{x}^{i''}\hat{y}^{p-1-i''}$.
Multiplying by $(\hat{x}+\hat{y})^{p-1}$ gives
\begin{equation}\label{E:hat}
  a_i(\hat{x}+\hat{y})^{p\lambda_i-1}=
  \sum_{i'=0}^{r-1}\sum_{j'=0}^{s-1}\alpha_{i,i',j'}\left(\sum_{i''=0}^{p-1}(-1)^{i''}\hat{x}^{i''}\hat{y}^{p-1-i''}\right)\hat{x}^{pi'}\hat{y}^{pj'}.
\end{equation}
However, the $prs$ basis elements $\hat{x}^{k}\hat{y}^{\ell}$ in~\eqref{E:hat}
are nonzero and distinct. Thus, since not all scalars $\alpha_{i,i',j'}$
are zero, it follows that $a_i(\hat{x}+\hat{y})^{p\lambda_i-1}\neq0$.
Hence $\dim(W_i)=p\lambda_i$ where
$W_i:=a_iF[\hat{x}+\hat{y}]$ for $1\leq i\leq r$.

We now prove that
$\hat{A}=\bigoplus_{i=1}^r\bigoplus_{j=0}^{p-1}(1+\hat{x})^jW_i$.
The following decompositions
\[
  F[\hat{x}]=\bigoplus_{i=0}^{p-1} \hat{x}^iF[\hat{x}^p]
  =\bigoplus_{i=0}^{p-1} (1+\hat{x})^iF[\hat{x}^p]\quad\textup{and}\quad
  F[\hat{y}]=\bigoplus_{j=0}^{p-1} \hat{y}^iF[\hat{y}^p]
  =\bigoplus_{j=0}^{p-1} (1+\hat{y})^jF[\hat{y}^p]
\]
imply that $\hat{A}=F[\hat{x},\hat{y}]\cong F[\hat{x}]\otimes_FF[\hat{y}]$
may be decomposed as
\[
  \hat{A}=\bigoplus_{i=0}^{p-1}\bigoplus_{j=0}^{p-1}(1+\hat{x})^i(1+\hat{y})^j
    F[\hat{x}^p]\otimes_FF[\hat{y}^p]
  =\bigoplus_{i=0}^{p-1}\bigoplus_{j=0}^{p-1}(1+\hat{x})^i(1+\hat{y})^jA.
\]
Now $(1+\hat{x})^{kp}(1+\hat{y})^{\ell p}A=(1+x)^k(1+y)^\ell A=A$ and so
$(1+\hat{x})^i(1+\hat{y})^jA=(1+\hat{x})^{i'}(1+\hat{y})^{j'}A$ holds
if $i\equiv i'\pmod p$ and  $j\equiv j'\pmod p$. Setting $k=i-j$ gives
\[
  \hat{A}=\bigoplus_{k=0}^{p-1}(1+\hat{x})^k\bigoplus_{j=0}^{p-1}
         (1+\hat{x})^j(1+\hat{y})^jA
         =\bigoplus_{k=0}^{p-1}(1+\hat{x})^k AF[(1+\hat{x})(1+\hat{y})],
\]
because $(1+\hat{x})^p(1+\hat{y})^p\in A$. However, we may replace
$F[(1+\hat{x})(1+\hat{y})]$ with $F[\hat{x}+\hat{y}]$ by~\eqref{E:pol}.
Using $A=\bigoplus_{i=1}^r a_iF[x+y]$ and
$F[x+y]\subseteq F[\hat{x}+\hat{y}]$ now gives
\[
  \hat{A}=\bigoplus_{k=0}^{p-1}(1+\hat{x})^k\bigoplus_{i=1}^ra_iF[x+y]F[\hat{x}+\hat{y}]
=\bigoplus_{i=1}^r\bigoplus_{k=0}^{p-1}a_i(1+\hat{x})^kF[\hat{x}+\hat{y}].
\]
Finally, $(1+\hat{x})^k$ is invertible, and so
$\dim(a_i(1+\hat{x})^kF[\hat{x}+\hat{y}])=p\lambda_i$, as desired.
\end{proof}

\section{Periodicity and duality}\label{S:PD}

Let $p^m$ be the smallest power of $p=\Char(F)$ satisfying $r\leq p^m$.
In this section we prove periodicity and a duality results
which depend on $p^m$.
The deviation vector $\eps(r,s,p):=(\eps_1,\dots,\eps_r)$ in Definition~1(c)
satisfies $\sum_{i=1}^r\eps_i=0$ since $\sum_{i=1}^r\lambda_i=rs$.
It turns out that periodicity and a duality are satisfied by
the deviation vector $\eps(r,s,p)$, but not the partition $\lambda(r,s,p)$.
The following lemma characterizes when $\eps(r,s,p)=(0,\dots,0)$;
Proposition~\ref{P:} follows from it.

\begin{lemma}\label{L:0}
Suppose that $\Char(F)=p$ and $r\leq\min\{p^m,s\}$.
\begin{itemize}
\item[(a)] If $s\equiv0\pmod{p^m}$, then $(x+y)^s=0$, $A=\bigoplus_{i=0}^{r-1} x^iF[x+y]$, and $x^i(x+y)^{s-1}\neq 0$ for $0\leq i<r$. Consequently,
$\eps(r,s,p)=(0,0,\dots,0)$.
\item[(b)] If $\eps(r,s,p)=(0,0,\dots,0)$, then $s\equiv0\pmod{p^m}$.
\end{itemize}
\end{lemma}

\begin{proof}
(a) Suppose that $s=kp^m$ where $k$ is an integer. Then
$(x+y)^{p^m}=x^{p^m}+y^{p^m}=y^{p^m}$ as $0=x^r=x^{p^m}$. Thus
$(x+y)^{kp^m}=y^{kp^m}=y^s=0$. It follows from $x^i(x+y)^s=0$ that
$\dim(x^iF[x+y])\leq s$. However,
$A=\sum_{i=0}^{r-1} x^iF[x+y]$ has dimension $rs$, and so the sum must be
direct. Thus $\dim(x^iF[x+y])=s$, and $x^i(x+y)^{s-1}\neq 0$
holds for $0\leq i<r$.

(b) It follows from $\eps_1=0$ that $\lambda_1=s$, and hence that $(x+y)^s=0$.
Equation~\eqref{E:bin} gives $\sum_{i=1}^{r-1} \binom{s}{i}x^iy^{s-i}=0$ and thus
$\binom{s}{i}=0$ in $\F_p$ for $1\leq i\leq r-1$.
A theorem of Lucas~\cite[p.\,2]{Gr} says that
$\binom{s}{i}\equiv\prod_{k\geq0}\binom{s_k}{i_k}\pmod p$
where $s=\sum_{k\geq0}s_kp^k$ and $i=\sum_{k\geq0}i_kp^k$ are the base-$p$ expansions
of $s$ and $i$, respectively. As $p^{m-1}<r\leq p^m$, we have $p^{m-1}\leq r-1$. Putting
$i=1,p,\dots,p^{m-1}$ into Lucas' theorem shows that $s_0=s_1=\cdots=s_{m-1}=0$.
In other words, $s\equiv0\pmod{p^m}$.
\end{proof}

The following lemma can be proved naturally using the theory of modules over
principal ideal rings, see~\cite{CK}. However, in the absence of a good
reference, our proof makes use of the more familiar theory of modules
over principal ideal domains.

\begin{lemma}\label{L:dual}
Suppose that $\overline{D}:=D/(\alpha^n)$  where $D$ is a principal ideal
domain and $\alpha\in D$ is prime. Suppose
$1\leqslant j\leqslant r$ and $M$ is a free $\overline{D}$-module
with basis $e_1,\dots,e_r$, and $N=\bigoplus_{i=1}^jx_i\overline{D}$
is a submodule of $M$ with
$\Ann(x_i)=(\alpha^{n_i})/(\alpha^n)\neq\overline{D}$ for $i=1,\dots,j$.
Then $M/N\cong\bigoplus_{i=1}^rD/(\alpha^{n-n_i})$ where $n_i=0$ for $j<i\leq r$.
\end{lemma}

\begin{proof}
We can (and will) initially view $M$ as a module over the
principal ideal domain~$D$.
Since $x_i\alpha^{n_i}=0$ for $i=1,\dots,j$, and $\alpha\in D$ is prime,
it follows from the theory
of $D$-modules \cite[Lemma\;9.1]{HH} that $x_i=y_i\alpha^{n-n_i}$
for some $y_i\in M$. Now view $M$ as a $\overline{D}$-module.
Let $(y_1,\dots,y_j)=(e_1,\dots,e_r)Y$
where $Y$ is an $r\times j$ matrix over $\overline{D}$. Then
\[
  \bigoplus_{i=1}^jy_i\overline{\alpha}^{\;n-1}\overline{D}
  =\bigoplus_{i=1}^j(y_i\overline{\alpha}^{\;n-n_i})\overline{\alpha}^{\;n_i-1}\overline{D}
  =\bigoplus_{i=1}^jx_i\overline{\alpha}^{\;n_i-1}\overline{D}
\]
is a $j$-dimensional linear space over the field
$\overline{D}/(\overline{\alpha})=D/(\alpha)$. We conclude
that there exists a $j\times j$ minor $Y_1$ of $Y$ such that $\det(Y_1)$ is
a unit in the local ring $\overline{D}$. (The set of units of
$\overline{D}$ equals $\overline{D}\setminus(\overline{\alpha})$
as $\alpha$ is prime in $D$.) Without loss of generality,
suppose that $Y_1$ comprises the first $j$ rows (and all $j$ columns) of $Y$.
Let $Y_2$ comprise the bottom $r-j$ rows of $Y$. Then
\begin{equation}\label{eq2}
  (y_1,\dots,y_j)=(e_1,\dots,e_j)Y_1+(e_{j+1},\dots,e_r)Y_2.
\end{equation}
Postmultiplying \eqref{eq2} by $Y_1^{-1}$ and rearranging gives
\begin{align*}
  (e_1,\dots,e_j)
    &=(y_1,\dots,y_j)Y_1^{-1}-(e_{j+1},\dots,e_r)Y_2Y_1^{-1}\textup{, and hence}\\
  (e_1,\dots,e_j,e_{j+1},\dots,e_r)&=(y_1,\dots,y_j,e_{j+1},\dots,e_r)
    \begin{pmatrix}Y_1^{-1}&0\\-Y_2Y_1^{-1}&I\end{pmatrix}.
\end{align*}
The above $r\times r$ matrix is invertible over $\overline{D}$.
Hence $y_1,\dots,y_j,e_{j+1},\dots,e_r$ is also a $\overline{D}$-basis of the
free $\overline{D}$-module $M$. With bases for $M$ and $N$ aligned,
it follows that
\[
  M/N\cong\bigoplus_{i=1}^rD/(\alpha^{n-n_i})
\]
where $n_i=0$ for $j<i\leq r$.
\end{proof}

We now prove $\eps(r,s,p)$ satisfies the periodicity and duality properties
in Theorem~\ref{T:DP}.

\begin{proof}[Proof of Theorem~\textup{\ref{T:DP}}]
Our strategy is to prove duality first, as duality implies periodicity.

(b) Suppose that $s=a+bp^m$ and $s'=-a+b'p^m$
where $a,b,b'$ are integers. Now $s'':=s+s'=(b+b')p^m$
is a multiple of $p^m$. Let $A$ be the homocyclic $F[x+y]$-module with
relations $x^r=y^{s''}=xy-yx=0$. By Lemma~\ref{L:0}, the nilpotent
transformation $\mu_{x+y}$ has minimal polynomial $t^{s''}$, and it corresponds
to the uniform Jordan partition $\lambda(r,s'',p)=(s'',\dots,s'')$.
The action of $\mu_{x+y}$ on $A$ gives
submodules of the same dimension as the action of $\mu_{(1+x)(1+y)}$ on $A$ as
discussed in Section~\ref{S:Basic}. The restriction of $\mu_{(1+x)(1+y)}$ to
the submodule $V_s$ of $A$ is $J_r\otimes J_{s''-s}=J_r\otimes J_{s'}$
by Table~\ref{T:Action}, and this
corresponds to the Jordan partition
$\lambda(r,s',p)=(s'+\eps'_1,\dots,s'+\eps'_r)$.
Similarly, the restriction of $\mu_{(1+x)(1+y)}$ to $A/V_s$ is $J_r\otimes J_s$
by Table~\ref{T:Action}, and this corresponds to the Jordan partition
$\lambda(r,s,p)=(s+\eps_1,\dots,s+\eps_r)$, say. Applying Lemma~\ref{L:dual}
with $M=A$, $N=V_s$, $n=s''$, $D=F[t]$, and $\alpha=t$ shows that
$\lambda(r,s,p)$ has parts $s''-(s'+\eps'_i)=s-\eps'_i$. Our ordering
conventions $\eps_1\geq\cdots\geq\eps_r$ and $\eps'_1\geq\cdots\geq\eps'_r$
imply that $\eps_{r-i+1}=-\eps'_i$. Thus $\eps(r,s,p)$ is
the negative reverse of $\eps(r,s',p)$. This proves part~(b).

(a) Suppose that $s\equiv s'\pmod{p^m}$. Choose an integer $s''$
such that $s''\geq r$ and $s''\equiv -s\equiv -s'\pmod{p^m}$. By part~(b),
$\eps(r,s'',p)$ is the negative reverse of both $\eps(r,s)$ and $\eps(r,s')$.
Hence $\overline{\eps(r,s,p)}=\overline{\eps(r,s',p)}$, and therefore
$\eps(r,s,p)=\eps(r,s',p)$.
\end{proof}

Henceforth, the phrase {\em by duality} will mean `by Theorem~\ref{T:DP}(b)',
and the phrase {\em by periodicity} will mean `by Theorem~\ref{T:DP}(a)'.

\begin{lemma}\label{L:bound}
If $r\leq s$ and $\eps(r,s,p)=(\eps_1,\dots,\eps_r)$, then
$|\eps_i|\leq r-1$ for $1\leq i\leq r$.
\end{lemma}

\begin{proof}
We noted in Section~\ref{S:Basic} that $(x+y)^{r+s-1}=0$.
Hence $\lambda_i\leq r+s-1$, and $\eps_i\leq r-1$ for $1\leq i\leq r$.
By duality, $-(r-1)\leq\eps_i$ and hence $|\eps_i|\leq r-1$ for $1\leq i\leq r$.
\end{proof}

Note that $|\eps_i|\leq\max\{|\eps_1|,|\eps_r|\}$
as $\eps_1\geq\cdots\geq\eps_r$. Proposition~\ref{P:s1} shows that the
upper bound of $r-1$ in Lemma~\ref{L:bound} can be attained.

\begin{proposition}\label{P:s1}
If $r\leq\min\{s,p^m\}$ and $s\equiv1\pmod{p^m}$, then
\[
  \eps(r,s,p)=(r-1,-1,\dots,-1).
\]
By duality, if $r\leq\min\{s,p^m\}$
and $s\equiv-1\pmod{p^m}$, then $\eps(r,s,p)=(1,\dots,1,-(r-1))$.
\end{proposition}

\begin{proof}
Suppose that $r\leq\min\{s,p^m\}$ and $s\equiv1\pmod{p^m}$. Let $q=p^{m+1}$. Since $\lambda(1,r,p)=(r)$, we have by \cite[(2.5a)]{G} that $\lambda(r,q-1,p)=(q,\dots,q,q-r)$. Subtracting the uniform vector $(q-1,q-1,...,q-1)$ gives $\eps(r,p^{m+1}-1,p)=(1,\dots,1,1-r)$. Then by duality, $\eps(r,s,p)=(r-1,-1,\dots,-1)$, as desired.
\end{proof}

\begin{proposition}\label{P:s2}
If $2\leq r\leq\min\{s,p^m\}$ and $s\equiv2\pmod{p^m}$, then
$$
\eps(r,s,p)=
\begin{cases}
(r-2,r-2,-2,\dots,-2)&\textup{if $r\equiv0\pmod{p}$,}\\
(r-1,r-3,-2,\dots,-2)&\textup{if $r\not\equiv0\pmod{p}$.}
\end{cases}
$$
By duality, $2\leq r\leq\min\{s,p^m\}$ and
$s\equiv-2\pmod{p^m}$ implies
$$
\eps(r,s,p)=
\begin{cases}
(2,\dots,2,2-r,2-r)&\textup{if $r\equiv0\pmod{p}$,}\\
(2,\dots,2,3-r,1-r)&\textup{if $r\not\equiv0\pmod{p}$.}
\end{cases}
$$
\end{proposition}

\begin{proof}
Suppose that $2\leq r\leq\min\{s,p^m\}$ and $s\equiv2\pmod{p^m}$.
Write $\eps(2,r,p)=(a,b)$. Then
$$
(a,b)=
\begin{cases}
(0,0)&\textup{if $r\equiv0\pmod{p}$,}\\
(1,-1)&\textup{if $r\not\equiv0\pmod{p}$,}
\end{cases}
$$
by Table \ref{T:eps}, and $\lambda(2,r,p)=(r+a,r+b)$.
Set $q=p^{m+1}$. 
By \cite[(2.5a)]{G}, we have $\lambda(r,q-2,p)=(q,\dots,q,q-r-b,q-r-a)$. Subtracting the uniform vector with all entries $q-2$ gives $\eps(r,q-2,p)=(2,\dots,2,2-r-b,2-r-a)$. Finally, duality shows that $\eps(r,s,p)=(r+a-2,r+b-2,-2,\dots,-2)$, and the proposition follows.
\end{proof}

\section{Which partitions \texorpdfstring{$\lambda(r,s,p)$}{} arise?}\label{S:small}

This section addresses the question: Which partitions $\lambda(r,s,p)$
arise?
Table~\ref{T:freq} lists the number, $n_r$, of deviation vectors
$\eps(r,s,p)$ as both $s\geq r$ and~$p$ vary.
Table~\ref{T:freq}, and parts of Table~\ref{T:eps}, were generated using
Magma computer code available at~\cite{Gl}. Theorem~\ref{T:nr} shows that
$n_r\leq 2^{r-1}$. Although this bound is optimal for $r\leq4$,
Table~\ref{T:freq} suggests that it may be a gross overestimate for large $r$.

\begin{table}[!ht]
\caption{The number, $n_r$, of different $\eps(r,s,p)$ vectors as $s\geq r$ and $p$ vary.}\label{T:freq}
\begin{tabular}{ |c|l|l|l|l|l|l|l|l|l|l|l|l| } \hline
$r$&1&2&3&4&5&6&7&8&9&10&11&12\\ \hline
$n_r$&1&2&4&8&14&24&28&45&61&78&94&118\\ \hline
\end{tabular}
\end{table}

The $i\times i$ matrix
$A_i=\begin{pmatrix}\binom{s+r-2i}{s-i+j-k}\end{pmatrix}_{0\leq j,k\leq i-1}$
in~\eqref{E:A} has integer entries.
Thus $\delta_i:=\det(A_i)$ is an integer, even though
the formulas~\eqref{E:delta2} and~\eqref{E:delta3} appear to give
rational values. As we are concerned with the case $\Char(F)=p>0$, we
henceforth assume that $\delta_i$, and the entries of $A_i$, lie
in the field~$\F_p$. The following recurrence for computing
$\lambda(r,s,p)$ is established in Section~2 of~\cite{II}.
%We believe the recurrence in Theorem~\ref{T:II} is correct, and have
%independently checked it gives correct values for $1\leq r\leq s\leq 20$ and $p\leq 11$.

\begin{theorem}[\cite{II}, Theorem 2.2.9]\label{T:II}
Suppose $r\leq s$ and $p=\Char(F)\geq0$. The parts
of $\lambda(r,s,p)$ can be computed recursively
\textup{(}in reverse order $\lambda_r,\lambda_{r-1},\dots,\lambda_1$\textup{)}
via
\begin{equation}\label{E:rec}
  \lambda_i=
  \begin{cases}
    r+s-2i+d(i)\quad&\textup{ if $\delta_i\neq0$,
     $\delta_{i-1}=\cdots=\delta_{i-(d(i)-1)}=0$ and $\delta_{i-d(i)}\neq0$},\\
    \lambda_{i+1}\quad&\textup{ if $\delta_i=0$}.
  \end{cases}
\end{equation}
\end{theorem}

When $\Char(F)=0$ the recurrence~\eqref{E:rec}
gives the familiar formula $\lambda_i=r+s+1-2i$, because
$d(i)=1$ and the sequence $\delta_{i-1},\dots,\delta_{i-(d(i)-1)}$ is empty.
Recall that $\delta_r=1$.

We now prove Theorem~\ref{T:nr}, which says for fixed $r\geq1$, that there are
at most $2^{r-1}$
different deviation vectors $\eps(r,s,p)$ as $s$, with $s\geq r$,
and the prime~$p$ vary.

\begin{proof}[Proof of Theorem~\textup{\ref{T:nr}}]
The recurrence relation in Theorem~\ref{T:II} for the $i$th part $\lambda_i$ of
$\lambda(r,s,p)$ depends on whether or not the determinants
$\delta_1,\dots,\delta_{r-1}$ are zero in $\F_p$. Subtracting~$s$
gives a recurrence relation for $\eps_i$ that depends on whether
$\delta_1,\dots,\delta_{r-1}$ are zero or nonzero, and is independent of $s$.
Hence there are at most $2^{r-1}$ choices for $\eps(r,s,p)$ as $s$,
with $s\geq r$,
and the prime~$p$ vary.
\end{proof}

Theorem~\textup{\ref{T:nr}} has computational implications for the construction
of Table~\ref{T:eps}. Theorem~\textup{\ref{T:finite}}, and it proof, gives
insight into the complexity of extending the $r$-values in Table~\ref{T:eps}.

\begin{proof}[Proof of Theorem~\textup{\ref{T:finite}}]
Fix $r$, and consider separately the cases: $p<r$, and $p\geq r$.

\textsc{Case $p<r$.}
By periodicity, we may assume that $s$ equals one of $r, r+1, \dots, r+p^m-1$.
This gives $p^m$ choices for $s$. Multiplying $p^{m-1}<r$ and $p<r$ gives
$p^m<r^2$. This case involves considering less than $r^2$ values of~$s$, and
less than $r$ values for the prime~$p$. Thus we must compute a finite number of
(less than $r^3$) deviation vectors $\eps(r,s,p)$ in this case.

\textsc{Case $p\geq r$.} By periodicity, we may assume that
$s\in\{r, r+1, \dots, r+p-1\}$ as $m=1$. Indeed, we may assume that
$s\in\{r, r+1, \dots, r+(r-2)\}$ by Theorem~\ref{T:std}, and by
duality.
For each of these $r-1$ choices for $s$, compute the $r-1$ integer determinants
$\delta_1,\dots,\delta_{r-1}$, and factor
$\Delta:=\delta_1\cdots\delta_{r-1}$ using~\eqref{E:delta3}.
Since the integers $\delta_1,\dots,\delta_{r-1}$ depend only on~$r$ and~$s$
and $r\leq s\leq 2r-2$, the number of primes $p$ dividing $\Delta$
is bounded by a function of~$r$. Thus the number of steps required
to compute $\eps(r,s,p)$ for $r\leq s\leq 2r-2$ and prime divisors~$p$
of $\Delta$ is bounded by a function of~$r$.
For the remaining primes $p$, the partition $\lambda(r,s,p)$ is standard, and
$\eps(r,s,p)=(r-1,r-3,\dots,-(r-3),-(r-1))$. Thus
a finite computation suffices to determine the values of $\eps(r,s,p)$
with $r\leq\min\{p,s\}$.
\end{proof}

We mention a subtle point: it is not a finite computation to determine a
mapping from the triples $(r,s,p)$
where $r$ is fixed and $s\geq r$ and $p$ vary, to the set of
of allowable values of $\eps(r,s,p)$. The latter requires the computation of
$s$ modulo~$p$ for infinitely many $s$ (and primes $p>r$), and this is an
infinite computation.

We conclude by stating an open problem.

\begin{problem}
Determine the asymptotic size as $r\to\infty$ of
the number $n_r$ of different vectors $\eps(r,s,p)$ where $s\geq r$
and~$p$ vary. (Does the $\lim_{r\to\infty}n_r/2^{r-1}$ exist?
{\it C.f.} Theorem~\textup{\ref{T:nr}}.)
\end{problem}

\noindent\textsc{Acknowledgements.}
We would like to thank Gary Seitz and Martin Liebeck for
helpful conversations. We thank Neil Strickland for his
answer to a question posed on MathOverflow\footnote{\href{http://mathoverflow.net/questions/134773/homocyclic-primary-module-over-pid}
{\tt http://mathoverflow.net/questions/134773/homocyclic-primary-module-over-pid}}
by the third author. We also thank Michael Barry for showing us his paper~\cite{B2}.
The first and
second authors acknowledge the support of the Australian Research
Council Discovery Grant DP110101153. This work was done during the visit of the third author to School of Mathematics and Statistics, University of Western Australia, and he would like to
thank the China Scholarship Council for its financial support.

\end{document}